\documentclass[12pt,a4paper]{amsart} 
\usepackage{enumerate,amsmath,amssymb}

\ifx\pdftexversion\undefined
\usepackage[dvips]{graphicx,color}
\else
\usepackage[pdftex]{graphicx,color}
\fi

\newcommand{\andname}{and}
\newcommand{\volumename}{volume}
\newcommand{\Inname}{in}
\newcommand{\ofname}{of}
\newcommand{\pagesname}{pages}
\newcommand{\deuxpoints}{:}
\newcommand{\EMdash}{---}

\let\ifanglais\iftrue

\def\R{{\mathbb R}}
\def\N{{\mathbb N}}

\def\vol{\mu}

\def\n#1{|| #1 ||}

\newcommand{\ed}[1]{\textrm{d} #1}

\newtheoremstyle{mesthm}
  {10pt plus 1pt minus 1pt}
  {9pt plus 1pt minus 6pt}
  {\slshape}
  {0.5cm}
  {\bfseries}
  {.}
  {1ex}
  {}
\newtheoremstyle{mesdefi}
  {6pt plus 1pt minus 1pt}
  {6pt plus 1pt minus 1pt}
  {}
  {0.5cm}
  {\bfseries}
  {.}
  {1ex}
  {}%

\theoremstyle{mesthm}
\newtheorem{lema}{\ifanglais{\large L}emma\else{\large L}emme\fi}
\newtheorem{theo}[lema]{\ifanglais{\large T}heorem\else {\large
    T}h\'eor\`eme\fi}

\newtheorem{prop}[lema]{{\large P}roposition}  
\newtheorem{cor}[lema]{{\large C}orollary}

\theoremstyle{mesdefi}

\newcounter{step}

\DeclareMathOperator{\ent}{Ent}

\title[Hilbert geometry of products]%
{On the Hilbert geometry of products}
\author[C. Vernicos]{Constantin Vernicos}
\address{%
  Institut de math\'ematique et de mod\'elisation de Montpellier\\ 
  Universit\'e Montpellier 2 \\
  Case Courrier 051\\
  Place Eug\`ene Bataillon \\
  F--34395 Montpellier Cedex\\ 
  France} 
\email{Constantin.Vernicos@math.univ-montp2.fr}

\begin{document}

\begin{abstract}
We prove that the Hilbert geometry of a product of convex sets is bi-lipschitz
equivalent the direct product of their respective Hilbert geometries.
We also prove that the volume entropy is additive with respect to product
and that amenability of a product is equivalent to the amenability of each terms.
\end{abstract}

\maketitle

\section*{Introduction and statement of results}

Hilbert geometries are simple metric geometries defined in the interior of a convex
set thanks to cross-ratios. They are generalisations of the projective
model of the Hyperbolic Geometry. Because of their definition they are invariant
by the action of projective transformations.
Among all these geometries, those admitting a discrete subgroup of their isometries acting
co-compactly, commonly known as \textsl{divisible Hilbert Geometries} or \textsl{divisible
convex sets}, play an
important part. For instance we can find examples of such geometries, which are Hyperbolic in
the sense of Gromov with a quotient which does not admit any Riemannian Hyperbolic metric \cite{benoist}. 

The present paper focuses on product of Hilbert Geometries, and takes its roots
in the following question: Does the product of two divisible convex sets give a divisible convex set ?
The answer to that question is no and is given by a very simple example, the Hilbert geometry of the square.
Indeed, the Hilbert geometry of the segment $[-1,1]$, which is isometric
to the real line, is divisible. However the product of two such segments, which is a square in $\R^2$, endowed
with its Hilbert geometry
is not a divisible convex set, which is related to the fact that one can't immerse $PGL(2,\R)\times PGL(2,\R)$
into $PGL(3,\R)$.
 
However following  B.~Colbois, C.~Vernicos, P.~Verovic~\cite{cvv2} the Hilbert
Geometry of a polygon is bi-lipschitz equivalent to $\R^2$.
In the light of that example we asked ourselves what is  the relation between the Hilbert Geometry of a product
and the product of Hilbert Geometries, our main theorem gives a complete answer to that question:

\begin{theo}[Main Theorem]
  The Hilbert Geometry of a product of open convex sets is bi-lipschitz equivalent to
the direct metric product of the Hilbert Geometries of those convex sets.
\end{theo}

The proof of that theorem is surprisingly simple but it allows us to get an impressive range
of corollaries. Noticeably with respect to the volume entropy (see also G.~Berck, A.~Bernig and C.~Vernicos~\cite{berck_Bernig_Vernicos}, and M.~Crampon~\cite{crampon})
and amenability (see C.~Vernicos~\cite{ver09}) we obtain the following consequences.

\begin{prop}[Main consequences]\label{maincorrolary}
  Consider the two bounded open convex sets $A$ and $B$, then
  \begin{enumerate}
   \item The volume entropy is additive : $$\ent(A\times B)= \ent(A)+\ent(B)\text{;}$$
   \item The product Hilbert geometry $(A\times B,d_{A\times B})$ is amenable if and only if both 
Hilbert geometries $(A,d_A)$ and $(B,d_B)$ are amenable.
  \end{enumerate}
\end{prop}

Although the product of divisible convex sets need not be divisible itself 
thanks to the Main consequences one can apply  M.~Crampon~\cite{crampon} theorem on volume entropy to get  new
rigidity results (see corollaries \ref{rigidityprod1} and \ref{rigidityprod2} in section \ref{entropy}).

Let us conclude by an "opening" remark. Our Main theorem  shows 
that a second family of Hilbert Geometry seems to play a similar
role to the divisible one, those we could call the \textsl{lip-divisible} ones, \textsl{i.e.}, whose
group of bi-Lipschitz bijections admits a discrete subgroup acting co-compactly. 
Noticeably, following our theorem, this family is closed under product.

\section{Definitions and notations}

A \textsl{proper} open set in $\R^n$ is a set not containing a whole line.

A Hilbert geometry
$(\mathcal{C},d_\mathcal{C})$ is a non empty \textsl{proper} open convex set $\mathcal{C}$
on $\R^n$ (that we shall call \textsl{convex domain}) with
the Hilbert distance 
$d_\mathcal{C}$ defined as follows: for any distinct points $p$ and $q$ in $\mathcal{C}$,
the line passing through $p$ and $q$ meets the boundary $\partial \mathcal{C}$ of $\mathcal{C}$
at two points $a$ and $b$, such that one walking on the line goes consecutively by $a$, $p$, $q$
$b$. Then we define
$$
d_{\mathcal C}(p,q) = \frac{1}{2} \ln [a,p,q,b],
$$
where $[a,p,q,b]$ is the cross ratio of $(a,p,q,b)$, i.e., 
$$
[a,p,q,b] = \frac{\| q-a \|}{\| p-a \|} \times \frac{\| p-b \|}{\| q-b\|} > 1,
$$
with $\| \cdot \|$ the canonical euclidean norm in
$\mathbb R^n$.  If either $a$ or $b$ is at infinity the corresponding ratio will be taken equal to $1$.

Note that the invariance of the cross ratio by a projective map implies the invariance 
of $d_{\mathcal C}$ by such a map.

These geometries are naturally endowed with
a  $C^0$ Finsler metric $F_\mathcal{C}$ as follows: 
if $p \in \mathcal C$ and $v \in T_{p}\mathcal C =\R^n$
with $v \neq 0$, the straight line passing by $p$ and directed by 
$v$ meets $\partial \mathcal C$ at two points $p_{\mathcal C}^{+}$ and
$p_{\mathcal C}^{-}$~. Then let $t^+$ and $t^-$ be two positive numbers such
that $p+t^+v=p_{\mathcal C}^{+}$ and $p-t^-v=p_{\mathcal C}^{+}$, in other words
these numbers corresponds to the time necessary to reach the boundary starting
at $p$ with the speed $v$ and $-v$. Then we define
$$
F_{\mathcal C}(p,v) = \frac{1}{2} \biggl(\frac{1}{t^+} + \frac{1}{t^-}\biggr) \quad \textrm{and} \quad F_{\mathcal C}(p , 0) = 0.
$$ 
Should $p_{\mathcal C}^{+}$ or
$p_{\mathcal C}^{-}$ be at infinity, then corresponding ratio will be taken equal to $0$.


The Hilbert distance $d_\mathcal{C}$ is the length distance associated to 
$F_{\mathcal C}$.

Thanks to that Finsler metric, we can built a Borel measure
$\mu_{\mathcal C}$ on $\mathcal C$ (which is actually
the Hausdorff measure of the metric space $(\mathcal C,
d_{\mathcal C})$, see \cite{bbi}, exemple~5.5.13 ) as follows.

To any $p \in \mathcal C$, let $B_{\mathcal C}(p) = \{v \in \R^n ~|~ F_{\mathcal{C}}(p,v) < 1 \}$
be the open unit ball in
$T_{p}\mathcal{C} = \R^n$ of the norm $F_{\mathcal{C}}(p,\cdot)$ and 
$\omega_{n}$ the euclidean volume of the open unit ball of the standard euclidean space
$\R^n$.
Consider the (density) function $h_{\mathcal C}\colon  \mathcal C \longrightarrow \R$ given by $h_{\mathcal C}(p)
= \omega_{n}/\text{Leb}\bigl(B_{\mathcal C}(p)\bigr),$ where $\text{Leb}$ is the canonical Lebesgue measure
of $\R^n$ equal to $1$ on the unit "hypercube". We define $\mu_{\mathcal C}$, which we shall call
the \textsl{Hilbert Measure} on $\mathcal{C}$,
by
$$
\mu_{\mathcal C}(A) = \int_{A} h_{\mathcal C}(p) \ed{\text{Leb}(p)}
$$
for any Borel set $A$ of $\mathcal C$.

The bottom of the spectrum of $\mathcal{C}$, denoted by $\lambda_1(\mathcal C)$,
and the Sobolev constant $S_\infty(\mathcal{C})$ are
defined as in a Riemannian manifold of infinite volume, thanks to
the Raleigh quotients as follows
\begin{equation}
\label{eqdef1}
\lambda_1(\mathcal{C})=
\inf \frac{\displaystyle{\int_\mathcal{C}{\n{df_p}_\mathcal{C}^*}^2~d\mu_\mathcal{C}(p)}}
{\displaystyle{\int_\mathcal{C} f^2(p)d\mu_\mathcal{C}(p)}}\text{,} \qquad%
 S_\infty(\mathcal{C}) = \inf \frac{\displaystyle{\int_\mathcal{C}{\n{df_p}_\mathcal{C}^*}~d\mu_\mathcal{C}(p)}}
{\displaystyle{\int_\mathcal{C} |f|(p)d\mu_\mathcal{C}(p)}}\text{,}
\end{equation}
where the infimum is taken over all non zero lipschitz functions with
compact support in $\mathcal{C}$

Finally the Cheeger constant of $\mathcal C$ is defined by
 \begin{equation}
\label{Cheeger}
I_{\infty}(\mathcal{C})=
\inf_U \frac{\nu_{\mathcal C}(\partial U)}{\mu_{\mathcal C}(U)}\text{,}
\end{equation}
where $U$ is an open set in $\mathcal{C}$ whose closure is compact and
whose boundary is a $n-1$ dimensional submanifold, and $\nu_{\mathcal C}$ is the
Hausdorff measure associated to the restriction of the Finsler norm $F_\mathcal{C}$
to hypersurfaces. Thanks to \cite{cvv0} we know that there is a constant $c$ such that
\begin{equation}
  \label{sobcheeger}
  \frac1c \cdot S_\infty(\mathcal{C})\leq I_\infty(\mathcal{C}) \leq c\cdot S_\infty(\mathcal{C})\text{.}
\end{equation}

\section{The decomposition lemma}\label{volume}

\begin{theo}\label{decomposition}
Consider the family of convex sets $A_i\in \R^{n_i}$, for $i=1,\ldots,k$ and $n_i\in \N^*$, then for any point $p=(p_1,\ldots,p_k)$ of the convex set $A_1\times\cdots\times A_k$ and 
any vector $v=(v_1,\ldots,v_k)\in \R^{n_1}\times \cdots \times \R^{n_k}$ one has
$$
\max_{1\leq i\leq k} F_{A_i}(p_i,v_i)\leq F_{A_1\times\cdots\times A_n}(p,v)\leq \sum_{i=1}^k F_{A_i}(p_i,v_i)\text{,}
$$
therefore the identity restricted to $A_1\times\cdots\times A_k$ 
is a bi-lipschitz map between $(A_1\times\cdots\times A_k,d_{A_1\times\cdots\times A_k})$ and the direct product of the metric spaces 
$(A_i,d_{A_i})$ for $i=1,\ldots,k$.
\end{theo}

\begin{proof}
  Consider a point  $p=(p_1,\ldots,p_k)$ of the convex set $A_1\times\cdots\times A_k$ and a vector  $v=(v_1,\ldots,v_k)\in \R^{n_1}\times \cdots \times \R^{n_k}$.
If the two positive numbers $t^+$ and $t^-$ are such that
$$
p+t^+v\in \partial(A_1\times\cdots\times A_k) \qquad \text{and} \qquad p-t^+v\in \partial(A_1\times\cdots\times A_k) 
$$
then $F_{A_1\times\cdots\times A_k}(p,v)=\frac12 \Bigl(\frac1{t^+} + \frac1{t^-}\Bigr)$.
This implies that for given $i,j\in\{1,\ldots,k \}$, $p_i+t^+v_i \in \partial A_i$ and $p_j-t^-v_j\in \partial A_j$.

Hence, should we define for each integer $i\in \{1,\ldots,k\}$ the positive numbers $t_i^+$ and $t_i^-$ by asking that
$$
p_i+t_i^+v_i\in \partial A_i \qquad \text{and} \qquad p_i-t_i^+v_i\in \partial A_i\text{,} 
$$
we would then obtain $t^+=\min\{t_1^+,\ldots,t_k^+\}$ and $t^-=\min\{t_1^-,\ldots,t_k^-\}$, which would imply that
$$
F_{A_1\times\cdots\times A_k}(p,v)=\frac12 \bigl(\max\{1/t_1^+,\ldots,1/t_k^+\} + \max \{1/t_1^-,\ldots,1/t_k^-\}\bigr) 
$$ 
and therefore using the classical comparison between the $l^1$ and $l^\infty$ norm in $\R^k$ we get
$$
\frac12 \max_{1\leq i\leq k} \biggl\{\frac1{t_i^+} + \frac1{t_i^-}\biggr\}\leq F_{A_1\times\cdots\times A_k}(p,v)\leq \frac12 \biggl(\sum_{i=1}^k \frac1{t_i^+} + \sum_{i=1}^k \frac1{t_i^-}\biggr)
$$
which we can rewrite by associativity of the addition in the following form:
$$
\max_{1\leq i\leq k} \frac12\biggl\{\frac1{t_i^+} + \frac1{t_i^-}\biggr\} \leq F_{A_1\times\cdots\times A_k}(p,v)\leq \sum_{i=1}^k \frac12 \Bigl(\frac1{t_i^+} +\frac1{t_i^-}\Bigr)\text{.}  
$$
\end{proof}

\begin{cor}\label{volumedecomposition}
  Consider the family of convex sets $A_i\in \R^{n_i}$, for $i=1,\ldots,k$ and $n_i\in \N^*$. Then at 
any point $p=(p_1,\ldots,p_k)\in A_1\times\cdots\times A_n$ we have the following inequality:
\begin{equation}
    \prod_{i=1}^k d\mu_{A_i}(p_i)\leq d\mu_{A_1\times\cdots\times A_n}(p)\\
\leq k^{n_1+\cdots+n_k}\prod_{i=1}^k d\mu_{A_i}(p_i)\text{.}
\end{equation}

\end{cor}

\begin{proof}
Let us denote by $\text{Leb}$ the Lebesgue measure on $\R^{n_1}\times\cdots\times\R^{n_k}$
normalised by $1$ on the unit cube 
and by $\text{Leb}_i$ the corresponding one on $R^{n_i}$.  Then, for any point $p=(p_1,\ldots,p_k)$ of the convex set $A_1\times\cdots\times A_k$,
if $T_pB(1)$ corresponds to the unit tangent ball at $p$ and for all $i$, $T_{p_i}B(1)$ to the unit tangent ball at $p_i$ for $A_i$ then one has
\begin{equation}
    \frac{1}{k^{n_1+\cdots+n_k}} \prod_{i=1}^k \text{Leb}_{i}\bigl(T_{p_i}B(1)\bigr)\leq\text{Leb}(T_pB(1))\\
\leq \prod_{i=1}^k \text{Leb}_{i}\bigl(T_{p_i}B(1)\bigr)
\end{equation}
and the corollary follows by definition of the measure.
\end{proof}

\section{Illustrations}
\begin{quote}
In order to illustrate in a simple way
  our main theorem we apply it to two geometries: the $n$-dimensional cube
and the $n$-dimensional simplex.
\end{quote}

The next two applications, are useful to obtain  qualitative information on volumes in the
given Hilbert geometries (see for instance proposition 6 in \cite{cvv1} and its corollaries 6.1 and 6.2).

\begin{prop}
  Let $\mathcal{C}^n=\mathopen]-1,1\mathclose[^n$ be the $n$-dimensional cube. We have the following
  \begin{enumerate}
  \item $\mathcal{C}^n$ is bi-lipschitz equivalent to $\R^n$.
  \item For all $x=(x_1,\ldots,x_n)\in \mathopen]-1,1\mathclose[^n$, let $T_xB(1)$ be the tangent unit ball for $F_{\mathcal{C}^n} $,
then we have
$$
(2/n)^n\prod_{i=1}^n(1-x_i^2)\leq\text{Leb}(T_xB(1))\leq 2^n\prod_{i=1}^n(1-x_i^2)\text{.}
$$
  \end{enumerate}
\end{prop}

Notice that actually there is a better lower bound because one can replace $(2/n)^n$ by $2^n/(n!)$, using theorem \ref{decomposition}
instead of its corollary.

\begin{prop}
  Let $\mathcal{S}^n=\mathopen]0,+\infty\mathclose[^n$ be the $n$-dimensional positive cone, 
whose Hilbert geometry is isometric to the Hilbert geometry of the simplex of $\R^n$. We have the following
\begin{enumerate}
\item $\mathcal{S}^n$ is bi-lipschitz equivalent to $\R^n$.
\item $x=(x_1,\ldots,x_n)\in \mathopen]0,+\infty\mathclose[^n$, let $T_xB(1)$ be the tangent unit ball for $F_{\mathcal{S}^n} $,
then we have
$$
(4/n)^n\prod_{i=1}^nx_i\leq\text{Leb}(T_xB(1))\leq 4^n\prod_{i=1}^nx_i\text{.}
$$
\end{enumerate}
\end{prop}

For the same reason, one can also replace $(4/n)^n$ by $4^n/(n!)$ in this lower bound.
In that case  one can actually  make a precise computation and obtain, for instance,
$$
\text{Leb}(T_xB(1))=12 x_1\cdot x_2
$$

\section{The volume entropy of products}\label{entropy}

\begin{quote}
  The general behaviour of the volume entropy is not yet completely understood,
and the main conjecture, to prove that it is always less than that of the Hyperbolic geometry, is still open in dimension bigger than $2$. Therefore the next result and the generalisation it implies validate this conjecture a little bit more. They also simplify and generalise result obtained in \cite{surent}
\end{quote}

\begin{prop}\label{additivityEnt}
The volume entropy is subadditive with respect to product of convex sets:
take a  family of bounded open convex sets $A_i\in \R^{n_i}$, for $i=1,\ldots,k$ and $n_i\in \N^*$, then one has
\begin{equation}
  \label{equnbounded}
  \max_{1\leq i\leq k}\bigl\{\ent(A_i)\bigr\}\leq\ent(A_1\times\cdots\times A_k)\leq \sum_{i=1}^k\ent(A_i)\text{.}
\end{equation}
If the convex sets are also bounded then we actually have additivity:
\begin{equation}
  \label{eqbounded}
  \ent(A_1\times\cdots\times A_k)=\sum_{i=1}^k\ent(A_i)\text{.}
\end{equation}
\end{prop}

\begin{proof}
  We will do the proof for $k=2$, the general case trivially follows.
Let $A$ and $C$ be two open convex sets, respectively in $\R^n$ and $\R^m$, and let $p=(p_A,p_C) \in A\times C$.

Thanks to the left hand side inequality of theorem \ref{decomposition}, we obtain
for any point $q=(q_A,q_c)\in  A\times C$ that
$$
d_A(p_A,q_A)\leq d_{A\times C}(p,q) \qquad \text{ and } \qquad d_C(p_C,q_C)\leq d_{A\times C}(p,q)\text{,}
$$
which imply the next inclusion
\begin{equation}
  \label{eqinclusionball}
  B_{A\times C}(p,R)\subset B_A(p_A,R)\times B_C(p_C,R)\text{.}
\end{equation}

The right hand side inequality of theorem \ref{decomposition} yields in turn that
for any $\varepsilon>0$ , $B_A(p_A,\varepsilon R)\times B_C\bigl(p_C,(1-\varepsilon)R\bigr)$ is a subset of $B_{A\times C}(p,R)$.


Hence, computing the volumes, using the inequalities of corollary \ref{volumedecomposition} we obtain
that
\begin{multline*}
 \vol_A\bigl(B_A(\varepsilon R)\bigr)\times \vol_C\bigl(B_C((1-\varepsilon)p,R)\bigr)\\
\leq  \vol_{A\times C}\bigl( B_{A\times C}(p,R)\bigr) \\
\leq 2^{n+m}\vol_A\bigl(B_A(p,R)\bigr)\times \vol_C\bigl(B_C(p,R)\bigr)
\end{multline*}

Taking the logarithm of both inequalities, dividing by $R$ and taking the limit as $R\to +\infty$, gives
the following inequality, for any $\varepsilon>0$:
$$
\varepsilon\ent(A)+(1-\varepsilon)\ent(C) \leq \ent(A\times C)\leq \ent(A) + \ent(C)\text{,}
$$
which implies \ref{equnbounded}.

In case both $A$ and $C$ are bounded, we can work as in \cite{cpv} with the \textsl{asymptotic balls} $\text{AsB}_{A\times C}(p,R)$, that is the
image of $A\times C$ by the dilation of ratio $\tanh(R)$ centred at $p$. Those asymptotic balls are exactly the
product of the asymptotic balls of $A$ and $C$ respectively centred at $p_A$ and $p_C$.
Therefore
\begin{equation}\label{eqcompact}
  \begin{split}
  \vol_{A}\bigl( \text{AsB}_{A}(p_A,R)\bigr)\times& \vol_{C}\bigl( \text{AsB}_{C}(p_C,R)\bigr)  \\
\leq &\vol_{A\times C}\bigl( \text{AsB}_{A\times C}(p,R)\bigr)\\ 
\leq &2^{n+m}\vol_{A}\bigl( \text{AsB}_{A}(p_A,R)\bigr)\vol_{C}\bigl( \text{AsB}_{C}(p_C,R)\bigr)\text{.}
\end{split}
\end{equation}

This inequality allows us to conclude using the fact shown in \cite{cpv} that there exists some constant $K$ such that
$$
B_{A\times C}(R-1)\subset\text{AsB}_{A\times C}(p,R)\subset B_{A\times C}(R+K)\text{.}
$$
\end{proof}

The following corollary is a straightforward application of M.~Crampon~\cite{crampon} rigidity result and the subadditivity of entropy:
\begin{cor}\label{rigidityprod1}
  Consider the family of divisible convex set with $C^1$ boundary $A_i\in \R^{n_i}$, for $i=1,\ldots,k$ and $n_i\in \N^*$, 
then one has
  \begin{itemize}
  \item $\ent(A_1\times\cdots\times A_k)\leq n_i-k$,
  \item Equality occurs if and only all $A_i$ are ellipsoids.
  \end{itemize}
\end{cor}

\begin{cor}
  Let $\mathcal{C}$ be a convex set in $\R^n\subset \R^{n+1}$ and let $p$ be a point outside $\R^n$ in $\R^{n+1}$.
Then  $\ent(p+\mathcal{C})=\ent(\mathcal{C})$.
\end{cor}
\begin{proof}
  This comes from the fact that $p+\mathcal{C}$ is projectively equivalent to $\mathcal{C}\times\mathopen]0,+\infty\mathclose[$,
and
$$
\max\bigl\{\ent C , \ent(\mathopen]0,+\infty\mathclose[)\bigr\}\leq \ent(\mathcal{C}\times\mathopen]0,+\infty\mathclose[)\leq \ent(\mathcal{C})+\ent(\mathopen]0,+\infty\mathclose[)$$ 
by Proposition \ref{additivityEnt}.
As $\mathopen]0,+\infty\mathclose[$ endowed with its Hilbert geometry is isometric to the real line we easily conclude.
\end{proof}

M.~Crampon~\cite{crampon} rigidity result applied to that case therefore implies:
\begin{cor}\label{rigidityprod2}
  Let $\mathcal{C}$ be a divisible convex set with $C^1$ boundary in $\R^n\subset \R^{n+1}$ and let $p$ be a point outside $\R^n$ in $\R^{n+1}$.
Then
\begin{itemize}
\item $\ent(p+\mathcal{C})\leq n-1$,
\item Equality occurs if and only if $\mathcal{C}$ is an ellipsoid.
\end{itemize}
\end{cor}

\section{Amenability of products}\label{amenability}

\begin{quote}
  Following our former work \cite{ver09} we say that a Hilbert geometry is amenable if and only if the bottom of its spectrum is null, which is equivalent to
the nullity of its Cheeger constant. In this section we show how this property behaves with respect to product.
\end{quote}

\begin{prop}
  Consider the family of convex sets $A_i\in \R^{n_i}$ for $i=1,\ldots,k$. The following are equivalent
  \begin{enumerate}[(i)]
    \item The Hilbert geometry of $A_1\times \cdots\times A_k$ is amenable;
    \item For all $i$, $A_i$ is amenable;
  \end{enumerate}
More precisely, with respect to the bottom of the spectrum and the Sobolev constants we have  the following inequalities:
\begin{eqnarray}
\label{eqinegalitelambda}  \lambda_1(A_1\times \cdots\times A_k)&\geq& k^{-n_1-\cdots-n_k}\max_{1\leq i\leq k} \lambda_1(A_i)\\
\label{eqinegalitecheeger} S_\infty(A_1\times\cdots\times A_k)&\leq& k^{n_1+\cdots+n_k} \sum_{1\leq i\leq k} S_\infty(A_i)\text{.}
\end{eqnarray}
\end{prop}
\begin{proof}
Let us denote $A_1\times \cdots\times A_k$ by $\Pi$.
  Consider a lipschitz function with compact support $f\colon A_1\times \cdots\times A_k\to \R$ , then we have for almost every point in $A_1\times \cdots\times A_k$ 
the function $f$ admits a differential $df$ and for any $i$ we have $||df||_{A_i}\leq ||df||_{\Pi}$, therefore for any $i$ we have
  \begin{equation}
    \label{eqinelambda}
    \int_{A_1\times \cdots\times A_k} ||df||^2_{\Pi} d\vol_{\Pi} \geq \int_{A_1\times \cdots\times A_k} ||df||^2_{A_i} d\vol_{\Pi}
  \end{equation}
Now thanks to Corollary \ref{volumedecomposition} we have for any $i$
\begin{eqnarray}
  \int_{\Pi} ||df||^2_{A_i} d\vol_{\Pi} &\geq& \int_{\Pi} ||df||^2_{A_i}  d\vol_{A_1}\cdots d\vol_{A_k}\\
  \text{then by definition of $\lambda$,}  &\geq & \lambda_1(A_i) \int_{\Pi} f^2  d\vol_{A_1}\cdots d\vol_{A_k}\text{,}\\
\text{and thanks to corollary \ref{volumedecomposition},}  &\geq & \dfrac{\lambda(A_i)}{ k^{n_1+\cdots+n_k}} \int_{\Pi} f^2 d\vol_{\Pi}\text{.}
\end{eqnarray}
which implies the inequality (\ref{eqinegalitelambda}), and the implication (i)$\Rightarrow$ (ii).

For the other implication we will use the Cheeger constant and for better clarity, restrict ourselves to the product
of two convex sets. 
Now let us suppose that $I(A)=I(C)=0$ and let us prove that $I(A\times C)=0$. To do so we will prove the inequality (\ref{eqinegalitecheeger}).
Let us
consider two real valued lipschitz functions $f$ and $g$ with compact support respectively in $A$ and $C$.
We then define the function $h\colon A\times C\to \R$ as follows: for any  $p=(p_A,p_C)\in A\times C$,   $h(p)=f(p_A)g(p_C)$. We first use
the textbook equality
$$
dh= gdf+fdg\text{.}
$$
Applying the right hand side inequality of Theorem  \ref{decomposition} we obtain
\begin{equation}
  ||dh||_{A\times C}\leq||dh||_A + ||dh||_C \leq |g|\cdot||df||_A +|f|\cdot||dg||_C\text{}.
\end{equation}

The next step consists in integrating over $A\times C$ taking into account the right hand side inequality of (\ref{volumedecomposition}) to obtain 
\begin{equation}
  \begin{split}
    \int_{A\times C} ||dh||_{A\times C} d\vol_{A\times C} &\leq \\2^{n_A+n_C}\biggl( &\int_C|g|d\vol_C \cdot\int_A||df||_Ad\vol_A
\\ &+ \int_A|f|d\vol_A \cdot\int_C||dg||_Cd\vol_C\biggr) \text.
  \end{split}
  \end{equation}

We finish by dividing by the integral of $|h|$ over $A\times C$ using the right hand side inequality of (\ref{volumedecomposition})
to finally get
\begin{equation}
  \begin{split}
    \dfrac{ \int_{A\times C} ||dh||_{A\times C} d\vol_{A\times C}}{\int_{A\times C} |h| d\vol_{A\times C}} \leq&\\ 2^{n_A+n_C}&\biggl(\dfrac{\int_A||df||_Ad\vol_A}{\int_A|f|d\vol_A} +\dfrac{\int_C||dg||_Cd\vol_C}{\int_C|g|d\vol_C}\biggr)\text{.}
  \end{split}
\end{equation}

This last inequality implies  inequality (\ref{eqinegalitecheeger}), and allows us to conclude thanks to 
the main theorem of our paper \cite{ver09}. 
\end{proof}

This proposition shades some light on the example given by proposition 4.1 in \cite{ccv}  of a Hilbert Geometry which is
not Hyperbolic in the sense of Gromov, but which has positive bottom of the spectrum,
and therefore allows us to get more example
of the same kind. Indeed, it is straightforward that a product of convex set is never strictly convex which implies that it is never 
Hyperbolic in the sense of Gromov.








\begin{thebibliography}{BBI01}




\bibitem[Ben06]{benoist}
\bgroup\bf Y.~Benoist\egroup{}.
\newblock {Convexes hyperboliques et quasiisom\'etries. (Hyperbolic convexes
  and quasiisometries.)}.
\newblock {\em Geom. Dedicata}, 122:109--134, 2006.

 \bibitem[BBI01]{bbi} \bgroup\bf D.~Burago\egroup{}, \bgroup\bf
   Y.~Burago\egroup{} \andname{} \bgroup\bf S.~Ivanov\egroup{}.
   \newblock \textsl{A Course in Metric Geometry}, \volumename{}~33
   \ofname{} \textsl{Graduate Studies in Mathematics}.  \newblock
   American Mathematical Society, 2001.
















\bibitem[BBV10]{berck_Bernig_Vernicos} \bgroup\bf G.~Berck\egroup{}, \bgroup\bf
  A.~Bernig\egroup{}  \andname{} \bgroup\bf
  C.~Vernicos\egroup{}.  \newblock Volume Entropy of Hilbert Geometries. 
\newblock \textsl{Pacific J. of Math.}  245(2):201--225, 2010.

\bibitem[Ber09]{andreas}
 \bgroup\bf A.~Bernig\egroup{}, \newblock {H}ilbert {G}eometry of {P}olytopes.
\newblock \textsl{Archiv der Mathematik}, 92:314-324, 2009.

\bibitem[CV06]{cvv0} \bgroup\bf B.~Colbois\egroup{} \andname{} \bgroup\bf C.~Vernicos\egroup{}.  
\newblock Bas du spectre et delta-hyperbolicit\'e en g\'eom\'etrie de Hilbert plane.
\newblock \textsl{Bulletin de la Soci\'et\'e Math\'ematique de France} 134(3):357--381, 2006.

\bibitem[CV07]{ccv} \EMdash   
\newblock Les g\'eom\'etries de Hilbert sont \`a g\'eom\'etrie locale born\'ee.
\newblock \textsl{Annales de l'institut Fourier} 57(4):1359-1375, 2007.

\bibitem[CV11]{cvv2} \bgroup\bf B.~Colbois\egroup{}, \bgroup\bf C.~Vernicos\egroup{}  
\andname{} \bgroup\bf P.~Verovic \egroup{}.  
 \newblock Hilbert Geometry for convex polygonal domains.
 \newblock preprint 2008, to appear in \textsl{Journal of Geometry}. arXiv:0804.1620v1 

 \bibitem[CVV04]{cvv1} \EMdash  \newblock L'aire des triangles id\'eaux en g\'eom\'etrie de Hilbert.
 \newblock \textsl{Enseign. Math.}, 50(3--4):203--237, 2004.


\bibitem[CV04]{cpv} \bgroup\bf B.~Colbois\egroup{} \andname{} \bgroup\bf
  P.~Verovic\egroup{}.  \newblock Hilbert geometry for strictly
  convex domain.  \newblock\textsl{Geom. Dedicata}, 105:29--42, 2004.

\bibitem[CV]{cvp2} \EMdash
 \newblock Hilbert domains quasi-isometric to normed vector spaces.
 \newblock preprint 2008,  arXiv:0804.1619v1 [math.MG].







\bibitem[Cra]{crampon} \bgroup\bf M.~Crampon\egroup{}. 
\newblock Entropies of compact strictly convex projective manifolds.
\newblock preprint arXiv:0904.2489, 2009.

\bibitem[dlH93]{dlharpe} \bgroup\bf P.~de~la Harpe\egroup{}.
  \newblock On {H}ilbert's metric for simplices.  \newblock \Inname{}
  \textsl{Geometric group theory, Vol.\ 1 (Sussex, 1991)},
  \pagesname{} 97--119. Cambridge Univ. Press, 1993.


 \bibitem[Hil71]{dhilbert} \bgroup\bf D.~Hilbert\egroup{}.  \newblock
   \textsl{Les fondements de la G\'eom\'etrie}, \'edition critique
     pr\'epar\'ee par P.~Rossier.  \newblock Dunod, 1971.


  



\bibitem[SM00]{so} \bgroup\bf {\'E}.~Soci{\'e}-M{\'e}thou\egroup{}.
  \newblock \textsl{Comportement asymptotiques et rigidit\'es en
    g\'eom\'etries de Hilbert},  \newblock th\`ese de doctorat de
  l'universit\'e de Strasbourg, 2000.  \newblock
  http\deuxpoints//www-irma.u-strasbg.fr/irma/pu\-bli\-ca\-tions/2000/00044.ps.gz.

\bibitem[SM02]{so2} \EMdash Caract\'erisation des ellipso\"\i des par
  leurs groupes d'automorphismes. \textsl{Ann. Sci. de l'{\'E}NS},
  35(4):537--548, 2002.

\bibitem[Ver09]{ver09} \bgroup\bf C.~Vernicos\egroup{}.
\newblock Spectral Radius and amenability in Hilbert Geometry.
 \newblock \textsl{Houston Journal of Math.}, 35(4):1143-1169, 2009.

\bibitem[Ver08]{surent}
   \EMdash
    \newblock Sur l'entropie volumique des g\'eom\'etries de Hilbert.
   \newblock {\em S\'em. Th. Spe. et Geo. de Grenoble}, 26:155-176, 2008.

\bibitem[Ver]{prepub08} \EMdash
\newblock Lipschitz Characterisation of Polytopal Hilbert Geometries.
\newblock arXiv:0812.1032v1 [math.DG], 2008.




\end{thebibliography}
\end{document}